\documentclass[a4paper]{amsart}
\usepackage{lineno,hyperref}
\usepackage{amsfonts}
\usepackage{amsthm}
\usepackage{amsmath}
\usepackage[all]{xy}
\usepackage{dsfont}
\usepackage{amssymb}
\usepackage{extarrows}
\usepackage{mathrsfs}
\usepackage{arydshln}

\usepackage{enumerate}

\newcommand{\C}{\ensuremath{\mathbb{C}}}

\newcommand{\Z}{\ensuremath{\mathbb{Z}}}

\newtheorem{theorem}{Theorem}[section]
\newtheorem{proposition}[theorem]{Proposition}
\newtheorem{corollary}[theorem]{Corollary}
\newtheorem{lemma}[theorem]{Lemma}

\theoremstyle{definition}
\newtheorem{definition}[theorem]{Definition}

\theoremstyle{remark}
\newtheorem{remark}[theorem]{Remark}
\newtheorem{example}[theorem]{Example}

\begin{document}
\title{(Co)Homology Self-closeness Numbers  of Simply-connected Spaces} 

\author{Pengcheng Li}          
\address{Academy of Mathematics and Systems Science, Chinese Academy of Sciences,
          University of Chinese Academy of Sciences}
   \email{xiaotianshishou@163.com} 
\date{}

\maketitle
\begin{abstract}
  The (co)homology self-closeness number of a simply-connected
  based CW-complex $X$ is the minimal number $k$ such that any
  self-map $f$ of $X$ inducing an automorphism of the (co)homology
  groups for dimensions $\leq  k$ is a self-homotopy equivalence. These
  two numbers are homotopy invariants and have a close relation
  with the group of self-homotopy equivalences. In this paper, we
  compare the (co)homology self-closeness numbers of spaces in
  certain cofibrations, define the mod $p$ (co)homology self-closeness
  number of simply-connected p-local spaces with finitely generated
  homologies and study some properties of the (mod $p$)
  (co)homology self-closeness numbers.

\end{abstract}

\section{Introduction}
The group of self-homotopy equivalences of a space  and its subgroups has been extensively studied by many mathematicians in history, such as  Arkowitz \cite{Ark2,Ark3}, Rutter \cite{Rut, Rut3}, Maruyama \cite{Maru}.  The groups of self-homotopy equivalences are usually difficult to compute.  In 2015 Choi and Lee \cite{HK} introduced the \emph{self-closeness number} $N\mathcal{E}(X)$ of a space $X$ to investigate the group $\mathcal{E}(X)$ of self-homotopy equivalences of $X.$ The  \emph{self-closeness number}  $N\mathcal{E}(X)$, which is denoted by  $N_\sharp\mathcal{E}(X)$ in this paper, is defined by 
\[ N_\sharp\mathcal{E}(X):=min\{k~|\mathcal{A}_\sharp^k(X)=\mathcal{E}(X)\},\]
where $\mathcal{A}_\sharp^k(X):=\{f\in [X,X]~|f_\sharp\colon \pi_i(X)\xrightarrow{\cong} \pi_i(X) \text{ for }i\leq k\}$.
Oda and Yamaguchi \cite{O.Y.} continued the study of self-closeness number and proved inequalities among the self-closeness numbers of spaces of a cofibration of the type:
\[S^{m+1}\xrightarrow{\gamma}B\xrightarrow{i}X\xrightarrow{p}S^{m+2}\] 
and gave  dual results of the comparison of self-closeness numbers of  spaces in a fibration of the type, \cite{O.Y.2}:
\[ K(G,m+1)\xrightarrow{q}X\xrightarrow{i}Y\xrightarrow{\gamma}K(G,m+2).\] 
Recently they published a paper involving the homology and cohomology self-closeness
numbers of a space, see Section 6 of \cite{OY3}. I avoid repeating the overlaps of some results and directly quote their results in this paper. The following notation is needed to make sense of the introduction.

We agree once for all that all spaces are simply-connected based CW-complexes and  maps are thought as the homotopy classes with the given representative. In notations, 
let $\mathcal{CW}_{sc}$ be the homotopy  category of  simply-connected based CW-complexes.  Let  $X, Y\in \mathcal{CW}_{sc}$,  $[X,Y]$  denotes the  set of homotopy classes of  based maps from $X$ to $Y$; identify a map $f$ with its homotopy classes ($f=[f]$) and $f=g$ means $f\simeq g$.  
Let  $H_i(X;G)$ be  the $i$-th reduced homology group of $X$ with  coefficient group $G$ and $H_i(X)=H_i(X;\Z)$. For a  map $f\colon X\to Y$, denote by  $f_\ast$ or $H_i(f;G)\colon H_i(X;G)\to  H_i(Y;G)$ the corresponding induced homomorphism. Similar notations are given for cohomologies. 

For simply-connected spaces,  the Whitehead theorem and the universal coefficient theorem for cohomology indicate that if  a map $f\colon X\to Y$ is a homotopy equivalence if and only if  $1.$  $f$ is a homology equivalence: $f_\ast\colon H_i(X)\xrightarrow{\cong} H_i(X)$ for  all $i$; or $2$ $f$ is a cohomology equivalence: $f^\ast\colon H^i(X)\xrightarrow{\cong} H^i(X)$ for all $i$. This motivates us to define the \emph{homology and cohomology self-closeness numbers}.

Let $X$ be a based CW-complex and consider the following subsets of $[X,X]$:
\begin{align*}
\mathcal{A}_\ast^k(X)&:= \{f\in [X,X]~| f_\ast\colon H_i(X)\xrightarrow{\cong} H_i(X)\text{~for } i\leq k\},\mathcal{A}_\ast^{\infty}(X):=\lim\limits_{k\to +\infty}\mathcal{A}_\ast^{k}(X);\\
\mathcal{A}^\ast_k(X)&:= \{f\in [X,X]~| f^\ast\colon H^i(X))\xrightarrow{\cong} H^i(X)\text{~for } i\leq k\}, \mathcal{A}^\ast_{\infty}(X):=\lim\limits_{k\to +\infty}\mathcal{A}^\ast_{k}(X).
\end{align*}
If $n\leq k$, by the Whitehead theorem  there is a chain of monoids by inclusion: 
\[\mathcal{E}(X)\subseteq\mathcal{A}_\ast^{\infty}(X)\subseteq \mathcal{A}_k^\ast(X)\subseteq \mathcal{A}^\ast_n(X)\subseteq  [X,X].\] There is a similar chain in the cohomology case.
  The \emph{homology self-closeness number} $N_\ast\mathcal{E}(X)$  and the \emph{cohomology self-closeness number } $N^\ast\mathcal{E}(X)$  of  $X$ are defined by:
\[N_\ast\mathcal{E}(X):=\min\{k~|\mathcal{A}_\ast^k(X)=\mathcal{E}(X)\}\text{ and } N^\ast\mathcal{E}(X):=\min\{k~|\mathcal{A}^\ast_k(X;\Z)=\mathcal{E}(X)\}.\]
They are both well-defined homotopy invariants (Proposition 37 of \cite{OY3}).

\begin{remark}\label{rmkdef}
\begin{enumerate}[$1$.]
\item If $X\in \mathcal{CW}_{sc}$, $N_\ast\mathcal{E}(X),N^\ast\mathcal{E}(X)$ take values in the range $\Z_{\geq 0}\cup \{+\infty\}$.  $\mathcal{E}_\ast(X)=0 $ if and only if $X$ is contractible, which is denoted by $X=\{\ast\}$.   $N^\ast\mathcal{E}(\bigvee_{n\geq 2}S^n)=N_\ast\mathcal{E}(\bigvee_{n\geq 2}S^n)=N_\sharp \mathcal{E}(\bigvee_{n\geq 2}S^n)=+\infty$ (Example 39 of \cite{OY3}).
\item  If $X$ is not simple or simply-connected, it  may happen that a self-map $f$ is a homology equivalence but not a homotopy equivalence, see Example 4.35 of \cite{AT}; in this case $\mathcal{A}_\ast^{k}(X)\neq \mathcal{E}(X)$ for any integer $k\geq 0$ and we denote by $N_\ast\mathcal{E}(X)= -\infty$. 

\end{enumerate}

\end{remark}

The connectivity degree of  $X$ is denoted by $conn(X)$, which means that $\pi_i(X)=0$ if $i\leq conn(X)$.  
  Let \[ H^\ast\text{-}dim(X):=\max\{i\geq 0~| H^i(X)\neq 0\},~H_\ast\text{-}dim(X):=\max\{i\geq 0~| H_i(X)\neq 0\}\] be the homology dimension and cohomology dimension of $X$, respectively.  It's easy to prove that if $\{\ast\}\neq X\in\mathcal{CW}_{sc}$, 
\[ conn(X)+1\leq N_\ast\mathcal{E}(X)\leq H_\ast\text{-}dim(X),~ conn(X)+1\leq N^\ast\mathcal{E}(X)\leq H^\ast\text{-}dim(X).\]

We can compare these three  types of self-closeness numbers of a simply-connected space and prove some inequalities among them, refer to Section 6 of \cite{OY3}. In this paper I choose the cohomology self-closeness number of simply-connected spaces to be the protagonist, since there are  richer structures in cohomology theory, such as the cohomology ring and the Steenrod operations. The paper is arranged as follows.

In Section 2, motivated  by Oda and Yamaguchi's paper \cite{O.Y.2}, I  quote some of Rutter's results  about extension of ladde rs of cofibrations  \cite{Rut2} and give a dual discussion on the cohomology self-closeness numbers of spaces in a generic cofibration $A\xrightarrow{\gamma}B\xrightarrow{i}X$.  By Theorem \ref{thmNCF1}, the inequality  $N^\ast\mathcal{E}(B)\leq N^\ast\mathcal{E}(X)$ holds  if the following conditions hold: 
\begin{enumerate}[$(1)$]
\item $n-1\leq conn(A)$,   $H^\ast\text{-}dim(B)\leq n-1$; 
\item  $\gamma\colon A\to B$ induces a surjection: $\gamma_\ast\colon [A,A]\to [A,B].$
\end{enumerate}
By Theorem \ref{thmNCF2}, $N^\ast\mathcal{E}(X)\leq N^\ast\mathcal{E}(B)$ holds under the following assumptions
\begin{enumerate}[$(1')$]
\item  $m-1\leq conn(B)<H^\ast$-$dim(B)\leq n-1\leq conn(A)< dim(A)\leq n+m-2$.
\item If there exists  maps $h\in [A, A]$ and  $g\in\mathcal{E}(B)$ such that $g\gamma=\gamma h$, then $h\in\mathcal{E}(A)$. 
\end{enumerate}
Moreover, if $(2')$ is substituted by the assumption that the induced map $\gamma_\ast\colon [A,A]\to [A,B]$ is bijective, then  $N^\ast\mathcal{E}(X)=N^\ast\mathcal{E}(B)$, see Theorem \ref{thmNCF3}. It follows that if $B$ is atomic ($N^\ast\mathcal{E}(B)=conn(B)+1$), then so is $X$(Corollary \ref{coratomic}). Consider the case $m=2, A=S^{n}$ and $m=3,A=P^{n+1}(q)=S^n\cup_qe^{n+1}$ respectively,  we get   Corollary \ref{corm=2} \ref{corm=3}; particularly, Corollary \ref{corm=2} is a cohomology version of Theorem 6 of \cite{OY3}. The special case where $[A,A]$ is a cyclic group $\Z$ or $\Z/q~(q\geq 2)$(Theorem \ref{thmspec}) can be viewed as a generalization Theorem 5 ($A=S^{n+1}$) of \cite{OY3}.

   In Section 3, we define the \textit{mod $p$ homology  self-closeness number }$N_\ast\mathcal{E}(X;p)$ and the \textit{mod $p$ cohomology  self-closeness number} $N^\ast\mathcal{E}(X;p))$ of a simply-connected $p$-local space $X$ with finitely generated homology. They are also well-defined homotopy invariants.
For such a space $X$, we have $N_\ast\mathcal{E}(X;p)=N^\ast\mathcal{E}(X;p)$(Proposition \ref{proph=ch}) and $N_\ast\mathcal{E}(X)=N_\ast\mathcal{E}(X;p)$(Proposition \ref{prop0p}). 

In Section 4 we prove some properties of (mod $p$) homology and cohomology self-closeness numbers.  Let $p$ be a prime or $p=0$, $X$ be a simply-connected  space with finitely generated homology if $p=0$ and further let $X$ be $p$-local if $p$ is a prime. Denote by $N^\ast\mathcal{E}(X;0)=N^\ast\mathcal{E}(X)$. By Proposition \ref{suspens}, \ref{propprod},  we have the following inequalities: \begin{align*}
&N^\ast\mathcal{E}(\Sigma X;p)\geq N^\ast\mathcal{E}(X;p)+1;\\
&N^\ast \mathcal{E}(X\times Y;p),~N^\ast \mathcal{E}(X\wedge Y;p), ~N^\ast \mathcal{E}(X\vee Y;p)\geq max\{N^\ast \mathcal{E}(X;p),N^\ast \mathcal{E}( Y;p)\}; \\
&N^\ast \mathcal{E}(\Sigma(X\times Y);p)\geq N^\ast \mathcal{E}(\Sigma(X\wedge	Y);p).
\end{align*}
The above inequalities are also true for (mod $p$) homology self-closeness numbers.
 If the cohomology ring $H^\ast(X;\Z/p)(\Z/0=\Z)$ is generated by the  classes $x_i\in H^{|x_i|}(X;\Z/p)$, $1,\cdots, m$,   then  by Proposition \ref{propchring} we have 
\[N^\ast\mathcal{E}(X;p)\leq max\{|x_1|,\cdots ,|x_m|\}.\]
Finally, I exhibit a result of Haibao Duan, Theorem \ref{duan}, which states  that for a simply-connected compact  K\"ahler manifold  $M$ with torsion free cohomology and $H^2(M)\cong\Z$, we have $N^\ast\mathcal{E}(M)=2$.

\section*{Acknowledgements}
The author would like to thank Professor Jianzhong Pan, and the editors Donald M. Davis and Martin Crossley, for revising some mistakes.

\section{Cohomology self-closeness number  and  cofibrations}
In this section we consider a generic cofibration $A\xrightarrow{\gamma}B\xrightarrow{i}X\xrightarrow{p}\Sigma A$ and discuss conditions for $N^\ast\mathcal{E}(B)\leq N^\ast\mathcal{E}(X)$ and $N^\ast\mathcal{E}(X)\leq N^\ast\mathcal{E}(B)$.
\subsection{Some lemmas}

\begin{lemma}\label{lemconn}
Let $f\colon X\to Y$ be a map between simply-connected spaces, $F_f$ be the homotopy fiber of $f$, $C_f$ be the cofiber of $f$. Then the following are equivalent:
\begin{enumerate}[$(i)$]
\item $f$ is $n$-connected.
\item $F_f$ is $(n-1)$-connected.
\item $C_f$ is $n$-connected.
\end{enumerate}
\begin{proof}
By Lemma 6.4.11 and Proposition 6.4.14 of \cite{Ark}.
\end{proof}
\end{lemma}

By the long exact sequence of (co)homology groups, the  five-lemma and the Whitehead theorem, it's clear that 
\begin{lemma}\label{lem3maps}
In the following  homotopy commutative diagram with fibering rows of  simply-connected spaces:
\[\xymatrix{
A\ar[d]^{h}\ar[r]^{\gamma}&B\ar[d]^{g}\ar[r]^{i}&X\ar[d]^{f}\\
A\ar[r]^{\gamma}&B\ar[r]^{i}&X
}\]
 if  any two vertical maps of $f,g,h$ are self-homotopy equivalences, so is the third one.
\end{lemma}

\begin{lemma}\label{lemext}
Let  $r\geq 1,n\geq 2$, $A\xrightarrow{\gamma} B\xrightarrow{i} X\xrightarrow{p}\Sigma A$ be a cofibration with $conn(A)\geq n-1$, $conn(X)\geq r$ and $dim(A)\leq r+n-1$. 
Given self-maps $g\colon B\to B$, $f\colon X\to X$ such that $fi=ig$, there exists a map $h\in $
such that   $\gamma h=g\gamma$.
\begin{proof}
A direct result of Proposition 4.4 of \cite{Rut2}.
\end{proof}
\end{lemma}

\begin{corollary}\label{corext}Let $r\geq n\geq 2$. If $n\leq conn(A)+1\leq dim(A)\leq r\leq conn(X)$, then given a self-map $g\colon B\to B$, there exist maps $f\colon X\to X$, $h\colon A\to A$ such that the following diagram is homotopy commutative:
\[\xymatrix{
A\ar[r]^{\gamma}\ar@{-->}[d]^{h}&B\ar[r]^{i}\ar[d]^{g}&X\ar@{-->}[d]^{f}\\
A\ar[r]^{\gamma}&B\ar[r]^{i}&X
}\]
\begin{proof}
The condition $dim(A)\leq conn(X)$ implies that the map 
\[i^\ast\colon [X,X]\to [B,X]\] is surjective, there exists a map $f\colon X\to X$ such that $fi=i^\ast(f)=i g$. Then apply Lemma \ref{lemext}.
\end{proof}
\end{corollary}

\begin{lemma}\label{lemext2}
Let  $m,n\geq 2$, $A\xrightarrow{\gamma} B\xrightarrow{i} X\xrightarrow{p}\Sigma A$ be a cofibration with $conn(A)\geq n-1, conn(B)\geq m-1$ and $dim(A)\leq m+n-2$.
Suppose there is a commutative diagram:
\[\xymatrix{
A\ar[r]^{\gamma}&B\ar[r]^{i}\ar[d]^{g}&X\ar[r]^{p}\ar[d]^{f}&\Sigma A\ar[d]^{\Sigma h}\\
A\ar[r]^{\gamma}&B\ar[r]^{i}&X\ar[r]^{p}&\Sigma A
}\]
Then $g\gamma=\gamma h$.
\begin{proof}
A direct result of Theorem 4.6 of \cite{Rut2}.
\end{proof}
\end{lemma}

\begin{corollary}\label{corext2}
Let $n\geq m\geq 2$, if $m\leq conn(B)+1\leq  dim(B)\leq n\leq conn(A)+1\leq dim(A)\leq n+m-2$, then given a map $f\colon X\to X$, there exists maps $h\colon A\to A$, $g\colon B\to B$ such that the following diagram is homotopy commutative, in which rows are cofibrations:
\[\xymatrix{
A\ar@{.>}[d]^{h}\ar[r]^{\gamma}&B\ar[r]^{i}\ar@{-->}[d]^{g}&X\ar[r]^{p}\ar[d]^{f}&\Sigma A\ar@{-->}[d]^{\Sigma h}\\
A\ar[r]^{\gamma}&B\ar[r]^{i}&X\ar[r]^{p}&\Sigma A
}\]
\begin{proof}
There condition $dim(A)\leq n+m-2\leq 2n-2\leq 2\cdot conn(A)$ implies that the suspension map 
\[\Sigma \colon [A,A]\to [\Sigma A,\Sigma A]\]
is bijective, by Theorem 1.21 of \cite{JC}. Then the result follows from Corollary \ref{corext} and Lemma \ref{lemext2}.
\end{proof}
\end{corollary}

Let $A\xrightarrow{\gamma}B\xrightarrow{i}X$ be a cofibration of simply-connected spaces. In the remainder of this section we shall investigate conditions for the comparison of  $N^\ast\mathcal{E}(B)$ and $N^\ast\mathcal{E}(X)$.

\subsection{Conditions for $N^\ast\mathcal{E}(B)\leq N^\ast\mathcal{E}(X)$}

\begin{theorem}\label{thmNCF1}
Let $n\geq 2$, $A\xrightarrow{\gamma}B\xrightarrow{i}X$ be a cofibration in $\mathcal{CW}_{sc}$. 
If the following conditions hold:
\begin{enumerate}[$(1)$]
\item $n-1\leq conn(A)$,   $H^\ast\text{-}dim(B)\leq n-1$.
\item $\gamma\colon A\to B$ induces a surjection: $\gamma_\ast\colon [A,A]\to [A,B]$.
\end{enumerate}
  Then $N^\ast\mathcal{E}(B)\leq N^\ast\mathcal{E}(X)$.
\begin{proof}
Since $N^\ast\mathcal{E}(B)\leq H^\ast\text{-}dim(B)\leq n-1$, we may suppose that $N^\ast\mathcal{E}(X)=k\leq n-1$ and $g\in\mathcal{A}^\ast_{k}(B)$.

 By the long exact sequence of cohomology groups and $conn(A)\geq n-1$, the induced homomorphism 
$i^\ast\colon H^{d}(X)\to H^d(B)$ is an isomorphism for $d\leq n-1$.

The surjectivity of $\gamma_\ast\colon [A,A]\to [A,B]$  implies that  there exists  a map $h\in [A,A]$ such that $g\gamma=\gamma h$ and hence there is a map $f\in [X,X]$ such that $fi=ig$.

Consider the following commutative diagram for $d\leq k\leq n-1:$
\[\xymatrix{
H^{d}(X)\ar[rr]^{i^\ast}\ar[rr]_{\cong}\ar[d]_{f^\ast}&&H^d(B)\ar[d]_{g_\ast}\ar[d]^{\cong}\\
H^{d}(X)\ar[rr]_{\cong}\ar[rr]^{i^\ast}&&H^d(B)
}\]
Then $g\in\mathcal{A}^\ast_k(B)$ implies that $f\in\mathcal{A}^\ast_k(X)=\mathcal{E}(X)$.

Since $H^\ast\text{-}dim(B)\leq n-1$,  the induced homomorphism 
\[\partial\colon H^{d-1}(A)\to H^d(X)\text{ is an isomorphism for $d\geq n+1$}.\]
Then the commutative diagram for $d\geq n+1:$
\[\xymatrix{
H^{d-1}( A)\ar[rr]^{\partial}\ar[rr]_{\cong}\ar[d]^{h^*}&&H^d(X)\ar[d]^{f^\ast}\ar[d]_{\cong}\\
H^{d-1}(A)\ar[rr]^{\partial}\ar[rr]_{\cong}&&H^d(X)
}\]
implies that $h^\ast\colon H^d(A)\to H^d(A)$ is an isomorphism for $d\geq n$ and hence for all $d\geq 0$, since $conn(A)\geq n-1$. By the Whitehead theorem, $h\in\mathcal{E}(A)$. Hence $g\in\mathcal{E}(B)$ by Lemma \ref{lem3maps} and therefore $N^\ast\mathcal{E}(B)\leq k=N^\ast\mathcal{E}(X)$.

\end{proof}
\end{theorem}

\subsection{Conditions for $N^\ast\mathcal{E}(X)\leq N^\ast\mathcal{E}(B)$}

Note that for a simply-connected CW-complex $B$ and $n \geq 2$, the condition
$H^\ast\text{-}dim(B) \leq  n -1$ implies $H_\ast\text{-} dim(B) \leq  n -1$, by the universal coefficient theorem
for cohomology. Then, by Proposition 4C.1 of \cite{AT}, $B$ admits a cell structure of
dimension at most $n$.
\begin{theorem}\label{thmNCF2}
Let $n, m\geq 2$,   $A\xrightarrow{ \gamma}B\xrightarrow{i}X$ be a cofibration. Consider the following assumptions: 
\begin{enumerate}[$(1)$]
\item  $m-1\leq conn(B)<H^\ast\text{-}dim(B)\leq n-1\leq conn(A)< dim(A)\leq n+m-2$.
\item If there exists  maps $h\in [A, A]$ and  $g\in\mathcal{E}(B)$ such that $g\gamma=\gamma h$, then $h\in\mathcal{E}(A)$. 
\end{enumerate}
If assumptions $(1)$ and $ (2)$ hold, then  $N^\ast\mathcal{E}(X)\leq N^\ast\mathcal{E}(B)$. 
\begin{proof}
Since $conn(A)\geq n-1$, $i\colon B\to X$ is $n$-connected.  The induced homomorphism \[i^\ast\colon H^d(X)\to H^d(B)\]
is an isomorphism for $d\leq n-1$ and an injection for $d=n$.

Suppose that $N^\ast\mathcal{E}(B)=k$ and $f\in\mathcal{A}^\ast_k(X)$. Then \[m\leq k\leq H^{\ast}\text{-}dim(B)\leq n-1.\]
 By Corollary \ref{corext2}, there exist self-maps $h\colon A\to A$, $g\colon B\to B$ filling in the  homotopy commutative diagram:
\[\xymatrix{
A\ar@{.>}[d]^{h}\ar[r]^{a\cdot \gamma}&B\ar[r]^{i}\ar@{-->}[d]^{g}&X\ar[r]^{p}\ar[d]^{f}&\Sigma A\ar@{-->}[d]^{\Sigma h}\\
A\ar[r]^{a\cdot \gamma}&B\ar[r]^{i}&X\ar[r]^{p}&\Sigma A
}\eqno(C1)\]
Consider the induced commutative diagram for $d\leq k\leq n-1:$
\[\xymatrix{
H^d(X)\ar[rr]^{i^\ast}\ar[rr]_{\cong}\ar[d]_{f^\ast}\ar[d]^{\cong}&&H^d(B)\ar[d]_{g^\ast}\\
H^d(X)\ar[rr]_{\cong}\ar[rr]^{i^\ast}&&H^d(B)
}\eqno(C2)\]
Then $f\in \mathcal{A}^\ast_{k}(X)$ implies that $g\in \mathcal{A}^\ast_k(B)=\mathcal{E}(B)$. By $(3)$ we then  have $h\in\mathcal{E}(A)$ and hence $f\in\mathcal{E}(X)$, by Lemma \ref{lem3maps}. Therefore,
 $N^\ast\mathcal{E}(X)\leq k=N^\ast\mathcal{E}(B)$.

\end{proof}

\end{theorem}

\begin{lemma}\label{lemEE}
If $\gamma\colon A\to B$ induces a bijection: 
$\gamma_\ast\colon [A,A]\to [A,B]$, then given a map $h\colon A\to A$ and  $g\in \mathcal{E}(X)$ such that $g\gamma=\gamma h$, we have $h\in\mathcal{E}(A)$.
\begin{proof}
  Let $\bar{g}\in\mathcal{E}$ be the homotopy inverse of $g$; that is, $\bar{g}g=1_B=g\bar{g}$. By the surjectivity of $\gamma_\ast \colon [A,A]\to [A,B]$, there exists a map $\bar{h}\in [A,A]$ satisfying $\bar{g}\gamma=\gamma \bar{h}$. We then have
  \[\gamma=\bar{g}g\gamma=\bar{g}\gamma h=\gamma\bar{h}h~~~~\text{ and }~~~~\gamma=g\bar{g}\gamma=g\gamma \bar{h}=\gamma h\bar{h}.\]
  Then, by the injectivity of $\gamma_\ast\colon [A,A]\to [A,B]$, we get $\bar{h}h=1_A=h\bar{h}$ and hence $h\in \mathcal{E}(B)$.
\end{proof}
\end{lemma}

\begin{theorem}\label{thmNCF3}
Let $n, m\geq 2$,   $A\xrightarrow{ \gamma}B\xrightarrow{i}X$ be a cofibration satisfying  the following conditions:
\begin{enumerate}[$(1)$] 
\item $m-1\leq conn(B)<H^\ast\text{-}dim(B)\leq n-1\leq conn(A)< dim(A)\leq n+m-2$.
\item $\gamma\colon A\to B$ induces a bijection: 
$\gamma_\ast\colon [A,A]\to [A,B]$.
\end{enumerate}
Then  $N^\ast\mathcal{E}(X)=N^\ast\mathcal{E}(B)$. 
\begin{proof}
$N^\ast\mathcal{E}(X)\geq N^\ast\mathcal{E}(B)$ by Theorem \ref{thmNCF1}; $N^\ast\mathcal{E}(X)\leq N^\ast\mathcal{E}(B)$ by Theorem \ref{thmNCF2} and Lemma  \ref{lemEE}.
\end{proof}
\end{theorem}

\begin{definition}
A CW-complex $X$ is called \emph{atomic }if $N_\sharp\mathcal{E}(X)=conn(X)+1$.

\end{definition}

It's immediate that 
\begin{lemma}
If $X\in\mathcal{CW}_{sc}$, the following are equivalent:
\begin{enumerate}[$(i)$]
\item $X$ is atomic.
\item $N^\ast\mathcal{E}(X)=conn(X)+1$.
\item $N_\ast\mathcal{E}(X)=conn(X)+1$.
\end{enumerate}
\end{lemma}

\begin{corollary}\label{coratomic}
Let $n, m\geq 2$,   $A\xrightarrow{ \gamma}B\xrightarrow{i}X$ be a cofibration satisfying  the following conditions:
\begin{enumerate}[$(1)$] 
\item $m-1\leq conn(B)<H^\ast\text{-}dim(B)\leq n-1\leq conn(A)< dim(A)\leq n+m-2$.
\item $\gamma\colon A\to B$ induces a bijection: 
$\gamma_\ast\colon [A,A]\to [A,B]$.
\end{enumerate}
If $B$ is atomic, then so is $X$.

\begin{proof}
Since $\Sigma A$ is $n$-connected, $i\ast\colon B\to X$ is $n$-connected, $i_\ast\colon\pi_i(B)\to \pi_i(X)$ is an isomorphism for $i\leq n-1$. Since $ m-1\leq conn(B)\leq n-1$, we have $conn(X)=conn(B)$ and hence
\[N^\ast\mathcal{E}(X)=N^\ast\mathcal{E}(B)=conn(B)+1=conn(X)+1.\]

\end{proof}
\end{corollary}

 Let $m=2$ and $A=S^n$, then we have
\begin{corollary}[A cohomological version of  Theorem 6 of \cite{OY3}]\label{corm=2}
Let $n\geq 2$, $a\neq 0$, $B$ be $1$-connected with  $H^\ast\text{-}dim(B)\leq n-1$ and let $S^{n}\xrightarrow{a\cdot \gamma}B\xrightarrow{i}X$ be a cofibration. If $\pi_{n}(B)\cong \Z\langle \gamma\rangle$, then $N^\ast\mathcal{E}(X)=N^\ast\mathcal{E}(B)$.

\end{corollary}

Let $m=3$ and $A=P^{n+1}(\Z/q)=S^n\cup_q e^{n+1}$, then we have
\begin{corollary}\label{corm=3}
Let $n\geq 3$,  $q\geq 2, (a,q)=1$, $B$ is $2$-connected with $H^\ast\text{-}dim(B)\leq n-1$ and let $P^{n+1}(\Z/q)\xrightarrow{a\cdot \gamma}B\xrightarrow{i}X$ be a cofibration. If $[P^{n+1}(\Z/q),B]\cong \Z/q\langle \gamma\rangle$, then $N^\ast\mathcal{E}(X)=N^\ast\mathcal{E}(B)$.

\end{corollary}

\begin{remark}
The above results are also true for homotopy and homology self-closeness after every $H^\ast\text{-}dim(B)$ is substituted by $H_\ast\text{-}dim(B)$ and every $N^\ast$ $N_\sharp$ and  $N_\ast$, respectively.
\end{remark}

\subsection{A special case}
Let $q\in \Z$. If $q>1$, denote the set of prime factors of $q$ by $Pr(q):$
\[Pr(q):=\{p_1,\cdots, p_l~| q=p_1^{r_1}\cdot \cdots\cdot p_l^{r_l}, p_i\text{ are primes },  r_i\geq 1\}.\]
\begin{theorem} \label{thmspec}
Let $n, m\geq 2$, $A=\Sigma^2 A'$ and let  $A\xrightarrow{a\cdot \gamma}B\xrightarrow{i}X\xrightarrow{p}\Sigma A$ be a cofibration with $a\cdot \gamma\in [A,B]$ nontrivial. If  the following assumptions hold: 
\begin{enumerate}
\item $m-1\leq conn(B)<H^\ast\text{-}dim(B)\leq n-1\leq conn(A)< dim(A)\leq n+m-2$.
\item $[A,A]\cong \Z/q\langle 1_A\rangle$ and $\gamma\in [A,B]$ is a generator of a direct summand $\Z/{q'}$, where $q,q'$ satisfy the conditions:
\[\left\{\begin{array}{ll}
q'=q,&q=0;\\[1ex]
q'|q, Pr(q')=Pr(q)&q\neq 0.
\end{array}\right.\eqno(C3)\]

\end{enumerate}
then  $N^\ast\mathcal{E}(X)\leq N^\ast\mathcal{E}(B);$ $``="$ holds if  further $[A,B]\cong \Z/q'\langle \gamma\rangle$.
\begin{proof}
For the inequality $N^\ast\mathcal{E}(X)\leq N^\ast\mathcal{E}(B)$, it suffices to show the new assumption $(2)$ above implies the ``old " $(2)$ in Theorem \ref{thmNCF2}. 

Let $\bar{g}$ be the inverse of $g$.  By the assumption $(2)$ we may put
 \[ h=s\cdot 1_A, \bar{g} \gamma=t\cdot \gamma+u\]
 for some $s\in\Z/q, t\in\Z/{q'}$ and $u\in [A,B]/{\Z/{q'}}$.
By $(C1)$ we have 
\[a\cdot\gamma=\bar{g} g (a\cdot \gamma)=ast\cdot \gamma+as\cdot u.\]
It follows that $s=1$ if $q=0$ and $s\equiv 1\pmod q$ if $q\neq 0$ by the condition $(C3)$. Thus $h\in \mathcal{E}(A)$ and therefore $f\in\mathcal{E}(X)$ by Lemma \ref{lem3maps}.
 
 If, in addition, $[A,B]\cong \Z/q'\langle \gamma\rangle$, we show that $N^\ast\mathcal{E}(B)\leq N^\ast\mathcal{E}(X)$.
 
 Suppose that $N^\ast\mathcal{E}(X)=l$ and $g\in\mathcal{A}^\ast_{l}(B)$. Then 
 \[ l=N^\ast\mathcal{E}(X)\leq N^\ast\mathcal{E}(B)\leq H^\ast\text{-}dim(B)=n-1.\]
 Let $h=s'\cdot 1_A\colon A\to A$. Since $[A,B]\cong \Z/q'\langle \gamma\rangle$, we have 
 $g\gamma=t'\cdot \gamma$ for some $t'\in\Z/{q'}$.
 Then there exist a map $f\colon X\to X$ such that
 \[ fi=ig, (\Sigma h)p=pf.\]
From the  commutative diagram $(C2)$ in the proof of Theorem \ref{thmNCF2} for $d\leq l$, we  see that $g\in\mathcal{A}^\ast_{l}(B)$ implies $f\in \mathcal{A}^\ast_{l}(X)=\mathcal{E}(X)$, which in turn implies $g\in\mathcal{E}(B)$ by the  commutative diagram $(C2)$ for $d\leq n-1$. Therefore $N^\ast\mathcal{E}(B)\leq N^\ast\mathcal{E}(X)$.

\end{proof}
\end{theorem}

Let $m=2$ and $A=S^n$. Then  we get a cohomological version of  Theorem 5 of \cite{OY3}.
\begin{corollary}\label{corThm3}
Let $n\geq 2$, $S^{n}\xrightarrow{a\cdot \gamma}B\xrightarrow{i}X$ be a cofibration, in which $B$ is $1$-connected and $H^\ast\text{-}dim(B)\leq n-1$. If  $0\neq a\in\Z, \gamma$ is a generator of a direct summand $\Z\subseteq\pi_{n}(B)$, then $N^\ast\mathcal{E}(X)\leq N^\ast\mathcal{E}(B)$.

\end{corollary}

Let $1_P$ be the identity of $P^{n}(q)$. It's well-known that if $q\equiv 1\pmod 2$, $n\geq 4$, \[ [P^{n}(q),P^{n}(q)]\cong \Z/q\langle 1_P\rangle.\]

   Let $m=3$ and $A=P^{n+1}(q)$, we have:
\begin{corollary}
Let $n\geq 2$,  $q,q'>1$ be odd integers such that $Pr(q)=Pr(q')$ and let $ P^{n+1}(q)\xrightarrow{a\cdot \gamma}B\xrightarrow{i}X$ be a cofibration, in which 
$B$ is $2$-connected and $H^\ast\text{-}dim(B)\leq n-1$.  If $\langle \gamma\rangle \subseteq \pi_{n+1}(B;\Z/q)\cong \Z/{q'}$ is a direct summand and  $a\cdot \gamma\neq 0$,   then $N^\ast\mathcal{E}(X)\leq N^\ast(B)$, and  equality holds if $\pi_{n+1}(B;\Z/q)\cong \Z/{q'}\langle \gamma\rangle.$

\end{corollary}

\subsection{Another condition for $N^\ast\mathcal{E}(B)\leq N^\ast\mathcal{E}(X)$ }

\begin{theorem}[a  cohomological version of Theorem 9 of \cite{O.Y.}.]\label{thmNCF0}
Let $r\geq n\geq 2$ and $A\xrightarrow{\gamma}B\xrightarrow{i}X$ be a cofibration. If one of the following conditions holds
\begin{itemize}
\item[$(1)$] $n\leq conn(A)+1\leq dim(A)\leq r\leq conn(X)$ and $ H^r(A)\cong H^r(B)$.
\item[$(2)$] $n\leq conn(A)+1\leq dim(A)<r\leq conn(X).$
\end{itemize}
then $N^\ast\mathcal{E}(B)\leq N^\ast\mathcal{E}(X)$.
\begin{proof}
Suppose that the condition $(2)$ holds. Since $X$ is $r$-connected and $H^r(A)=0$, by the long exact sequence of cohomology groups,  the induced homomorphism 
\[\gamma^\ast\colon H^d(B)\to H^d(A)\text{ is an isomorphism for } d\leq r.\eqno(C4)\]
Since $dim(A)<r$, the induced homomorphism
 \[i^\ast\colon H^d(X)\to H^d(B)\text{ is an isomorphism for } d\geq r+1.\eqno(C5)\]

If $(1)$ holds, we can also get the above $(C4), (C5)$.

Suppose that $N^\ast\mathcal{E}(X)=k\geq conn(X)+1\geq r+1$ and $g\in \mathcal{A}^\ast_k(B)$. 
By Corollary \ref{corext},  there exist self-maps $f\colon X\to X$ and $h\colon A\to A$ such that \[fi=ig, g\gamma=\gamma h.\]

Consider the following commutative diagram:
\[\xymatrix{
H^{d}(X)\ar[rr]^{i^\ast}\ar[d]^{f^\ast}&&H^d(B)\ar[d]^{\cong}\ar[rr]^{\gamma^\ast}\ar[d]_{g^\ast}&&H^d(A)\ar[d]_{h^\ast}\\
H^{d}(X)\ar[rr]^{i^\ast}&&H^d(B)\ar[rr]^{\gamma^\ast}&&H^d(A)
}\]
Then for $d\leq r<k$, by the second square above and $(C4)$, $g\in \mathcal{A}^\ast_k(B)\subseteq \mathcal{A}^\ast_r(B)$ implies that $h\in \mathcal{A}^\ast_r(A)=\mathcal{E}(A)$.
For $r+1\leq d\leq k$, by the first square above and  $(C5)$, $g\in \mathcal{A}^\ast_k(B)$ implies that $f^\ast\colon H^d(X)\to H^d(X)$  is an isomorphism. Since $conn(X)\geq  r$, we get $f\in \mathcal{A}_k^\ast(X)=\mathcal{E}(X)$. 
Hence $g\in \mathcal{E}(B)$  and therefore $N^\ast\mathcal{E}(B)\leq k=N^\ast\mathcal{E}(X)$.

\end{proof}

\end{theorem}

\section{mod $p$  (co)homology self-closeness numbers}

Let  $p$ be a prime, $\Z/p$ be the set of  integers modulo $p$, $\Z_{p}$ be the set of integers localized at $p$.  Let $\mathcal{CW}_{scpft}$ be the category of simply connected $p$-local CW-complexes with finitely generated homology group over $\Z_{p}$ in each dimension. 

 We shall use the following  universal coefficient theorem for cohomology: 
\begin{lemma}\label{lemUCT}
 For each $i\geq 1$ and  a CW-complex $X$, there is an isomorphism:
\[ H^i(X;\Z/p)\xrightarrow{\cong}\mathrm{Hom}_{\Z/p}(H_i(X;\Z/p),\Z/p).\]

\end{lemma}

There is an easier criterion to determine a homotopy equivalence in $\mathcal{CW}_{scpft}:$ 
\begin{lemma}\label{lemhch}
Let $p$ be a prime and  $f\colon X\to Y$ be a map (morphism) in the category $\mathcal{CW}_{scpft}$. Then the following are equivalent:
\begin{itemize}
\item[$(1)$] $f$ is a homotopy equivalence.
\item[$(2)$] $f_\ast \colon H_i(X;\Z/p)\to H_i(Y;\Z/p)$ is an isomorphism for all $i\geq 0$.
\item[$(3)$] $f^\ast\colon H^i(Y;\Z/{p})\to H^i(X;\Z/{p})$ is an isomorphism for all $i\geq 0$.
\end{itemize}
\begin{proof}
$(1)\iff (2)$ is a restatement of Lemma 1.3 of \cite{Wilkerson};  $(2)\iff (3)$ by Lemma \ref{lemUCT}.

\end{proof}
\end{lemma}
Hence for $X\in\mathcal{CW}_{scpft}$,  we can detect  self-homotopy equivalences of $X$ by the induced automorphisms of  $H_i(X;\Z/p)$ or $H^i(X;\Z/p)$.
\begin{definition}\label{defNp}
Let  $X\in \mathcal{CW}_{scpft}$. 
\[\mathcal{A}^\ast_{k}(X;p):=\{f\in [X,X]~|f^\ast\colon H^i(X;\Z/p)\xrightarrow{\cong} H^i(X;\Z/p) \text{ for } i\leq k\}.\]
The \emph{mod-$p$  cohomology self-closeness number}   $N^\ast\mathcal{E}(X;p)$ is defined by :
\[ N^\ast\mathcal{E}(X;p):=min\{k~|\mathcal{A}^\ast_k(X;p)=\mathcal{E}(X)\}.
\]
The monoids $\mathcal{A}_\ast^{k}(X;p) $ and the $\emph{mod-$p$  homology self-closeness number }   $ $N_\ast\mathcal{E}(X;p)$ are defined after replacing  cohomology by homology. 
\end{definition}

  It's easy to see that  $N^\ast\mathcal{E}(X;p), N^\ast\mathcal{E}(X;p)$ are  homotopy invariants,  by a parallel proof of Proposition 37 of \cite{OY3}. 
  
  \begin{proposition}\label{proph=ch}
Let $p$ be a prime,  $X\in\mathcal{CW}_{scpft}$. Then $N^\ast\mathcal{E}(X;p)=N_\ast\mathcal{E}(X;p)$.
 \begin{proof}
 By Lemma \ref{lemUCT} we have $\mathcal{A}^\ast_k(X;p)=\mathcal{A}_\ast^k(X;p)$ for each $k\geq 0$. Then the equality in the proposition follows.
 
 \end{proof} 
 \end{proposition}

\begin{proposition}\label{propp=0}
Let $X\in \mathcal{CW}_{sc}$ such that $H_i(X)$ is finitely generated for each $i$. Then $N_\ast\mathcal{E}(X)\leq  N^\ast\mathcal{E}(X)\leq N_\ast\mathcal{E}(X)+1;$ $N_\ast\mathcal{E}(X)=N^\ast\mathcal{E}(X)$ if $H^{k+1}(X)$ is free for $k=N_\ast\mathcal{E}(X)$.

\begin{proof}
By Proposition  41, 43, 45 of \cite{OY3}.
\end{proof}
\end{proposition}

\begin{example}\label{exmoore}
Let $q\neq 0, n\geq 2$, then  
$N_\ast\mathcal{E}(P^n(q))=n-1<N^\ast\mathcal{E}(P^n(q))=n$.
\end{example}

\begin{proposition}\label{prop0p}
 Let $p$ be  a prime and $X\in\mathcal{CW}_{scpft}$,  then  $N_\ast\mathcal{E}(X)=N_\ast\mathcal{E}(X;p)$.

\begin{proof}
Suppose that  $f\colon X\to X$. By the naturality of the universal coefficient theorem for homology, there is a commutative diagram for each $k:$
\[\xymatrix{
H_k(X)\otimes\Z/p~\ar@{>->}[r]\ar[d]^{f_\ast\otimes \Z/p}&H_k(X;\Z/p)\ar@{->>}[r]\ar[d]^{f_\ast}&\mathrm{Tor}^{\Z}_1(H_{k-1}(X),\Z/p)\ar[d]^{\mathrm{Tor}(f_\ast,\Z/p)}\\
H_k(X)\otimes\Z/p~\ar@{>->}[r]&H_k(X;\Z/p)\ar@{->>}[r]&\mathrm{Tor}^{\Z}_1(H_{k-1}(X),\Z/p)
}\]
It follows that $\mathcal{A}_\ast^k(X)\subseteq \mathcal{A}_\ast^k(X;p)$ and hence $N_\ast\mathcal{E}(X)\leq N_\ast\mathcal{E}(X;p)$.

Suppose that $N_\ast\mathcal{E}(X)=l$ and $f\in\mathcal{A}_\ast^l(X;p)$. Then by the long exact sequence of homology groups we have $H_i(C_f;\Z/p)=0$ for $i\leq l$ and hence $H_i(C_f)\otimes \Z/p=0$ for $i\leq l$, by the universal coefficient theorem for homology. Since $X$ is $p$-local, so is  $C_f$. It follows that $H_i(C_f)=0$ for $i\leq l$ and hence the homomorphism $f_\ast\colon H_i(X)\to H_i(X)$ is an isomorphism for $i\leq l-1$ and an epimorphism for $i=l$. Since $H_k(X)$ is finitely generated, $f\in\mathcal{A}_\ast^l(X)=\mathcal{E}(X)$. Therefore $N_\ast\mathcal{E}(X;p)\leq l=N_\ast\mathcal{E}(X)$.

\end{proof}
\end{proposition}

\begin{example}\label{exChang}
Let $n\geq 3,t,r\geq 1$, $C^{n+2,t}_r=P^{n+1}(2^r)\cup_{i\eta q} \mathbf{C}P^{n+1}(2^t)$ be the four-cell \text{Chang complex}, where $\eta\in\pi_1^s$ is the suspension of the Hopf map, $i\colon S^n\to P^{n+1}(2^t)$ and $q\colon P^{n+1}(2^r)\to S^{n+1}$ are the canonical inclusion and quotient maps. We have
\[N_\ast\mathcal{E}(C^{n+2,t}_r)=N_\ast\mathcal{E}(C^{n+2,t}_r;2)=n.\]
\begin{proof}
The proof is parallel  to that of Lemma 3.1 of \cite{Z.L.}.
\end{proof} 

\end{example}

\section{More properties of self-closeness numbers}
  Let $\mathcal{CW}_{sc0ft}$ be the category of simply connected  CW-complexes with finitely generated homology group in each dimension. $\Z/0=\Z$.
We temporarily adopt the following \textbf{notations:}
\begin{align*}
 & \mathcal{A}^\ast_k(X;0):=\mathcal{A}^\ast_k(X),N^\ast\mathcal{E}(X;0):=N^\ast\mathcal{E}(X); \\
 &\mathcal{A}_\ast^k(X;0):=\mathcal{A}_\ast^k(X),N_\ast\mathcal{E}(X;0):=N_\ast\mathcal{E}(X).
\end{align*}

\begin{proposition} \label{suspens}
Let $p$ be a prime or $p=0$ and $\{\ast\}\neq X\in  \mathcal{CW}_{scpft}$. Then 
\begin{enumerate}[$(1)$]
\item $N^\ast\mathcal{E}(\Sigma X;p)\geq N^\ast\mathcal{E}(X;p)+1;$ equality holds if $dim(X)\leq 2\cdot conn(X)+1$.
\item $N_\ast\mathcal{E}(\Sigma X;p)\geq N_\ast\mathcal{E}(X;p)+1;$ equality holds if $dim(X)\leq 2\cdot conn(X)+1$. 
\end{enumerate}

\begin{proof} 
$(1)$ Suppose that  $N_\ast\mathcal{E}(\Sigma X;p)=k+1$ for some $k\geq 0$ and   $f\in\mathcal{A}_\ast^k(X;p)$.  By the natural isomorphism $H^{i}(X;\Z/p)\xrightarrow{}H^{i+1}(\Sigma X;\Z/p)$,
 $\Sigma f\in \mathcal{A}_\ast^{k+1}(\Sigma X;p)=\mathcal{E}(\Sigma X)$.    By the naturality again, we get 
\[f^\ast\colon H^i(X;\Z/p)\xrightarrow{\cong} H^i(X;\Z/p),~\forall i\geq 0.\]
Thus $f\in\mathcal{E}(X)$ by Lemma \ref{lemhch} and therefore $N^\ast\mathcal{E}(X;p)\leq k=N^\ast\mathcal{E}(\Sigma X;p)-1$.

If $dim(X)\leq 2\cdot conn(X)+1$, by Theorem 1.21 of \cite{JC}, the suspension map
\[ \Sigma\colon [X,X]\xrightarrow{} [\Sigma X, \Sigma X]\] is  a surjection.
Suppose that $N^\ast\mathcal{E}(X;p)=l$ and $F\in \mathcal{A}^\ast_{l+1}(\Sigma X;p)$ such that $F=\Sigma f$ for some $f\in [X,X]$. Then we have $f\in \mathcal{A}^\ast_l(X;p)=\mathcal{E}(X)$ and hence $F=\Sigma f\in \mathcal{E}(\Sigma X)$. Thus $N^\ast\mathcal{E}(\Sigma X;p)\leq l+1=N^\ast\mathcal{E}(X;p)+1$.

$(2)$ The proof of $(2)$  is completed after replacing ``cohomology" with ``homology" in $(1)$ above.

\end{proof}

\end{proposition}

 It's  easy to get  $N_\sharp\mathcal{E}(\C \emph{P}^n)=N_\ast\mathcal{E}(\C \emph{P}^n)=N^\ast\mathcal{E}(\C \emph{P}^n)=2$.
\begin{example}\label{excpn} 
$N_\sharp\mathcal{E}(\Sigma\C \emph{P}^2)=N_\ast\mathcal{E}(\Sigma\C \emph{P}^2)=N^\ast\mathcal{E}(\Sigma \C \emph{P}^2)=5$.
\begin{proof} Write $C^{5}_\eta=\Sigma\C \emph{P}^2=S^3\cup_{\eta}e^5$.
 By Theorem 41, 45 of \cite{OY3}, 
 we have \[N_\sharp\mathcal{E}(C^{5}_\eta)=N^\ast\mathcal{E}(C^{5}_\eta)=N_\ast\mathcal{E}(C^{5}_\eta).\]

By Section 8 of  \cite{Toda}, $[C^{5}_\eta,C^{5}_\eta]\cong \Z\langle 1_\eta\rangle\oplus\Z\langle i_3\bar{\zeta}\rangle$, where $1_\eta$ is the identity of $C^{5}_\eta$, $i_3\colon S^3\to C^{5}_\eta$ is the canonical inclusion map and  $ \bar{\zeta}\in [C^{5}_\eta,S^3]$ and $\tilde{\zeta}\in [S^{5},C^{5}_\eta]$ satisfy the relations (relations (8.3) and (8.4) of \cite{Toda}):
\[ \bar{\zeta}i_3=2\cdot 1_{3}, ~~~~q_{5}\tilde{\zeta}=2\cdot 1_{5},~~~~ i_3\bar{\zeta}+\tilde{\zeta}q_{5}=2\cdot 1_\eta,\eqno(C6)\]
where $1_n$ is the identity of $S^n$ and $q_5\colon C^{5}_\eta\to S^5$ is the canonical quotient map. 

 Let $\sigma_n1$ be the image of $1\in H_0(S^0)$ under the suspension: $H_0(S^0)\xrightarrow[\cong]{\Sigma^n}H_n(S^n)$.  We have  
 \[H_k(C^5_\eta)\cong \left\{\begin{array}{cl}
  \Z~\langle a_\eta\rangle,& k=3;\\
  \Z~\langle b_\eta\rangle,&k=5;\\
  0,&otherwise.
\end{array} \right.\] 
where $a_\eta= (i_3)_\ast(\sigma_3 1), b_\eta=(q_5)_\ast^{-1}(\sigma_5 1)$.
 It follows that  $N_\ast\mathcal{E}(C^5_\eta)=3$ or $5$.

By the relations $(C6)$, it's easy to get that \[(\overline{\zeta})_\ast(a_\eta)=2\cdot \sigma_31,~~~~ (i_3\overline{\zeta})_\ast(b_\eta)=0;~~~~(\widetilde{\zeta})_\ast(\sigma_{5} 1)=2\cdot b_\eta.\]
 We compute that  $f=x\cdot 1_\eta+y\cdot i_n\overline{\zeta}\in \mathcal{A}_\ast^{3}(C^{5}_\eta)$ with $x,y\in \Z$   if and only if $x+2y=\pm 1$. Note that $f=3\cdot 1_\eta-i_3\overline{\zeta}\notin \mathcal{E}(C^{5}_\eta):$ $f_\ast(b_\eta)=3\cdot b_\eta$, so  we get 
 \[\mathcal{E}(C^{5}_\eta)=\mathcal{A}_\ast^{5}(C^{5}_\eta)\subsetneqq  \mathcal{A}_\ast^3(C^{5}_\eta), ~~~~N_\ast\mathcal{E}(C^{5}_\eta)=5.\]
 
\end{proof}
\end{example}

\begin{proposition}\label{propprod}
Let $p$ be a prime or $p=0$ and $X,Y\in\mathcal{CW}_{scpft}$. Then
\begin{enumerate}[$(1)$]
\item $N^\ast\mathcal{E}(X\vee Y;p)\geq max\{N^\ast\mathcal{E}(X;p),N^\ast\mathcal{E}(Y;p)\}$.
\item $N^\ast\mathcal{E}(X\wedge Y;p), N^\ast\mathcal{E}(X\times Y;p)\geq max\{N^\ast\mathcal{E}(X;p),N^\ast\mathcal{E}(Y;p)\}$.
\item $N^\ast\mathcal{E}(\Sigma (X\times Y);p)\geq N^\ast\mathcal{E}(\Sigma (X\wedge Y);p)$.
\end{enumerate} 
Similar results hold for mod $p$ homology self-closeness numbers.
\begin{proof}
$(1)$ Assume that $N^\ast\mathcal{E}(X\vee Y)=k<max\{N^\ast\mathcal{E}(X),N^\ast\mathcal{E}(Y)\}=N^\ast\mathcal{E}(X)$. Suppose  $f\in \mathcal{A}_k^\ast(X;p), g\in \mathcal{A}^\ast_k(Y;p)$.  By the  natural isomorphism: $H^d(X\vee Y;\Z/p)\cong H^d(X;\Z/p)\oplus H^d(Y;\Z/p)$, we have $f\vee g\in \mathcal{A}^\ast_{k}( X\vee Y;p)=\mathcal{E}(X\vee Y).$
It follows that $f\in \mathcal{E}(X),g\in\mathcal{E}(Y)$ and hence $\mathcal{A}^\ast_k(X)=\mathcal{E}(X)$. Therefore $ N^\ast\mathcal{E}(X)\leq k= N^\ast\mathcal{E}(X\vee Y)$, which contradicts with the assumption.

$(2)$ The proof is similar  to that  of Proposition 46 $(2)$ of \cite{OY3}, using the general K\"unneth formula for cohomology with coefficients $\Z/p$.

$(3)$ By Proposition 4I.1 of \cite{AT}, there is a homotopy equivalence: 
\[\Sigma (X\times Y)\simeq \Sigma X\vee \Sigma Y\vee \Sigma (X\wedge 	Y).\]
Then the inequality follows from $(1)$ and $(2)$.

Replacing $N^\ast$ by $N_\ast$ and ``cohomology" by the dual  ``homology", we get  the proof  the corresponding results of mod $p$ homology self-closeness numbers.
\end{proof}
\end{proposition}

\begin{example}\label{exsphere}  Let $n\geq m\geq 2$.  We have 
\[N^\ast\mathcal{E}(\Sigma (S^m\times S^n))=m+n+1>N^\ast\mathcal{E}(S^m\times S^n)=N^\ast\mathcal{E}(S^m\vee S^n)=n.\]

\begin{proof} By Proposition \ref{propprod} we have
\[m+n+1\geq N^\ast\mathcal{E}(\Sigma (S^m\times S^n))\geq N^\ast\mathcal{E}(\Sigma (S^m\wedge  S^n))=m+n+1.\]

$N^\ast\mathcal{E}(S^m\times S^n)=N^\ast\mathcal{E}(S^m\vee S^n)=n$ follows from Theorem 41, Theorem 45 of \cite{OY3} and Proposition 5 of \cite{O.Y.}. One can also prove this by  applying Theorem \ref{thmNCF3} to the cofibration:
 \[ S^{m+n-1}\xrightarrow{[i_1,i_2]} S^m\vee S^n\xrightarrow{i}S^{m}\times S^{n},\]
 where $i_1\colon S^m\hookrightarrow S^m\vee S^n$, $i_2\colon S^n\hookrightarrow S^m\times S^n$ are the canonical inclusion maps and $[i_1,i_2]\in \pi_{m+n-1}(S^m\vee S^n)$ is their Whitehead product, a generator of a direct summand $\Z$.

\end{proof}
\end{example}

\begin{proposition}\label{proplocal}
 Let $p$ be a prime or $p=0$,  $X\in \mathcal{CW}_{sc}$ and  $l_p\colon X\to X_p$ be the localization at $p$, then
 \begin{enumerate}[$(1)$]
 \item $N_\ast\mathcal{E}(X)\leq max\big\{N_\ast\mathcal{E}(X_p)~|p \in\{\emph{primes},0\}\big\}\leq H_\ast\text{-}dim(X).$
 \item $N_\sharp\mathcal{E}(X)\leq max\big\{N_\sharp\mathcal{E}(X_p)~|p \in\{\emph{primes},0\}\big\}\leq H_\ast\text{-}dim(X)+1.$ \\
 If, in addition, $H_n(X)$ is finitely generated for $n=H_\ast\text{-}dim(X)$, then \[ max\big\{N_\sharp\mathcal{E}(X_p)~|p \in\{\emph{primes},0\}\big\}\leq H_\ast\text{-}dim(X).\]
 \item If $X$ is a torsion space $(X_0=\{\ast\}\text{ or }(\pi_i(X)\otimes \mathbb{Q}=0)$, then 
 \[N_\square\mathcal{E}(X)=max\big\{N_\square\mathcal{E}(X_p)~|p \in\{\emph{primes}\}, \square=\ast,\sharp.\]
\end{enumerate}  
\begin{proof}
$(1)$ Suppose that $max\big\{N_\ast\mathcal{E}(X_p)~|p \in\{\text{primes},0\}\big\}=k$ and $f\in\mathcal{A}_\ast^k(X)$. For each  $ p\in\{\text{primes},0\}$, by the universal property of localization, there is a unique (up to homotopy) map $f_p\colon X_p\to X_p$ such that $l_p f=f_pl_p$. Consider the following commutative diagram:
\[\xymatrix{
H_i(X)\ar@/^+6mm/[rr]^{-\otimes \Z_{p}}\ar[r]^{l_{p\ast}}\ar[d]^{f_\ast}&H_{i}(X_{p})\ar[r]^{\cong~~~~}\ar[d]^{f_{p\ast}}&H_i(X)\otimes \Z_{p}\ar[d]^{f_\ast\otimes \Z_{p}}\\
H_i(X)\ar@/_+6mm/[rr]_{-\otimes \Z_{p}}\ar[r]^{l_{p\ast}}&H_{i}(X_{p})\ar[r]^{\cong~~~~}&H_i(X)\otimes \Z_{p}
}\]
Since $-\otimes\Z_{p}$ is an exact functor,  $f\in\mathcal{A}_\ast^k(X)$ implies that $f_\ast\otimes\Z_{p}$ is an isomorphism for $i\leq k$ and hence $f_{p}\in \mathcal{A}_\ast^k(X_{p})=\mathcal{E}(X_{p})$ for all $p$.
Thus $f\in \mathcal{A}_\ast^{\infty}(X)=\mathcal{E}(X)$ and $N_\ast\mathcal{E}(X)\leq k= max\big\{N_\ast\mathcal{E}(X_p)~|p \in\{\text{primes},0\}\big\}$.
 
 For the second $``\leq"$, since $H_i(X_p)\cong H_i(X)\otimes\Z_{p}$ for each prime $p$ or $p=0$, we have
 \[N_\ast\mathcal{E}(X_p)\leq H_\ast\text{-}dim(X_p)\leq H_\ast\text{-}dim(X).\]

$(2)$ The proof of the first $``\leq "$ is similar and the second $``\leq "$ follows from Theorem 3 of \cite{OY3}.

$(3)$ If $X$ is a torsion space, by $(5)$ of \cite{loc} (page 41), there is a product decomposition:
\[X=X_{p}\times	 \prod\limits_{q\neq p} X_{q}.\]
Thus  
$N_\square\mathcal{E}(X)\geq max\big\{N_\square\mathcal{E}(X_p)~|p \in\{\text{primes}\}\big\}$, by Proposition \ref{propprod} if $\square=\ast$ and  by  Theorem 3 of \cite{HK} if $\square=\sharp$.
\end{proof}
\end{proposition}

\begin{proposition}\label{propchring}
Let $p$ be a prime or $p=0$, $X\in\mathcal{CW}_{scpft}$.
If the cohomology ring $H^\ast(X;\Z/p)$ is generated by the  cohomology classes $x_i\in H^{k_i}(X;\Z/p)(i=1,\cdots,m)$ with $k_1\leq \cdots\leq k_m$, then 
$N^\ast\mathcal{E}(X;p)\leq k_m$. 
\begin{proof}
 Suppose that $f\in \mathcal{A}^\ast_{k_m}(X;p)$. Then the induced ring homomorphism \[f^\ast\colon H^\ast(X;\Z/p)\to H^\ast(X;\Z/p)\] is surjective, since all generators $x_i$ are in the image.  
 Then each grading $H^i(X;\Z/p)$ is finitely generated implies that the induced epimorphism $f^\ast\colon  H^i(X;\Z/p)\to H^i(X;\Z/p)$ is an isomorphism for all $i$.  Thus $f\in\mathcal{A}^\ast_{\infty}(X;p)=\mathcal{E}(X)$ by Lemma \ref{lemhch}.

\end{proof}
\end{proposition}

\begin{lemma}\label{lemmfd}
Let $M$ be a closed simply-connected  manifold of dimension $2n$.  If $f\colon M\to M$ is a map  of degree $\pm 1$, then $f\in\mathcal{E}(M)$ if and only if $f\in \mathcal{A}^\ast_n(M)$.
\begin{proof}
Suppose that $f\in\mathcal{A}^\ast_n(M)$.
By 12 Theorem(p.248) of \cite{Spanier}, there is a natural short exact sequence:
\[ 0\to \mathrm{Ext}(H^{i+1}(M),\Z)\rightarrow H_i(M)\rightarrow Hom(H^i(M),\Z)\to 0.\]
Hence $f\in\mathcal{A}^\ast_n(M)$ implies that $f\in \mathcal{A}_\ast^{n-1}(M)$.
Then by the natural  Poinc\'{a}re duality $H^{n+i}(M)\cong H_{n-i}(M)$ for $i=1,\cdots,n$, we get $f\in\mathcal{A}^\ast_{2n-1}(M)$. Since $deg(f)=\pm 1$, $f^\ast\colon H^{2n}(M)\to H^{2n}(M)$ is an isomorphism and hence $f\in\mathcal{A}^\ast_{2n}(M)=\mathcal{E}(M)$.

\end{proof}
 \end{lemma}
 
 We end the paper with a theorem given by Professor Haibao Duan.

\begin{lemma}[The Hard  Lefschetz Theorem]\label{lef}
 Let $M$ be a simply-connected compact K\"ahler manifold of real dimension $2n$ with the K\"ahler class $d\in H^2(M,\mathbb{Q})$.  
Then the multiplication
\[ d^{n-r}\cup-\colon H^r(M;\mathbb{Q})\rightarrow H^{2n-r}(M;\mathbb{Q})\]
is an isomorphism for $0\leq r\leq n$.
\end{lemma}

\begin{theorem}[Duan]\label{duan}
Let $M$  be a simply-connected compact K\"ahler manifold with torsion free cohomology and $H^2(M)$ a cyclic group, then $N^\ast\mathcal{E}(M)=2$.
\begin{proof}Let  $dim(M)=2n$.  We may choose a K\"ahler class $d$ of $M$  such that $(M,d)$ is a K\"ahler manifold with $H^2(M;\Z)\cong \Z\langle d\rangle$. 
 By  Lemma \ref{lemmfd}, it suffices to show that a self-map $f$ of $M$ satisfying $f^\ast(d)=\varepsilon\cdot  d(\varepsilon=\pm 1)$ belongs to $\mathcal{A}^*_n(M)$.
 
 For each $2\leq r\leq n$, since $H^r(M)$ is torsion  free,  there exist cohomology classes $x_1,\cdots,x_{m_r}$ such that  $H^r(M)\cong \bigoplus_{i=1}^{m_r}\Z\langle x_i\rangle$. Then $\{x_i\}_{i=1}^{m_r}$ is also a $\mathbb{Q}$-basis of $H^r(M;\mathbb{Q})$. By  Lemma \ref{lef}, $\{d^{n-r}x_i\}_{i=1}^{m_r}$ is  a basis of $H^{2n-r}(M;\mathbb{Q})$. 
 There are relations:
 \[d^{n-r}x_ix_j=a_{ij}d^n, a_{ij}\in \mathbb{Q},1\leq i,j\leq m_r.\]
 Then $A=(a_{ij})_{m_r\times m_r}$ is a non-singular matrix by the  Poinc\'{a}re duality.
 
 Let $f^\ast(x_i)=\sum_{k=1}^{m_r}b_{ik}x_k$ and put $ B_r=(b_{ij})\in M_{m_r}(\Z)$.  Applying the ring homomorphism $f^\ast$ to the above relations, we have \[
 \varepsilon^{n-r}d^{n-r}(\sum_{k=1}^{m_r}b_{ik}x_k)(\sum_{k=1}^{m_r}b_{jk}x_k)=a_{ij}\varepsilon^nd^n.\]
Let $B_r^T$ denotes the transpose of $B_r$, we get an equality of matrices:
 \[\varepsilon^{n-r}B_rAB^{T}_r=\varepsilon^nA.\]
The non-singularity of $A$ then implies that $det(B_r)^2=\varepsilon^{rm_r}=1$. Thus $B_r$ is non-singular  and  therefore $f\in \mathcal{A}^*_n(M)$.
 
\end{proof}
\end{theorem}

\begin{example}\label{exGrass}
Let $ n\leq m<\infty$,  $G_n(\mathbb{C}^m)$ be the Grassmannian of $n$-dimension vector subspaces of $\mathbb{C}^m.$ By 4.10 Example of \cite{Kahler}, $G_n(\mathbb{C}^m)$ is a K\"ahler manifold and by  Chapter 6, 14 of \cite{Milnor},  $G_n(\mathbb{C}^m)$ satisfies other conditions in Theorem \ref{duan}. Thus
$N^\ast(G_n(\mathbb{C}^m))=2$.

\end{example}

\bibliographystyle{siam}
\bibliography{references}
\end{document}